\documentclass[a4paper,12pt]{amsart}

\usepackage[utf8]{inputenc}
\usepackage{amsmath,amssymb,amsthm,amscd,amsfonts}
\usepackage{latexsym}

\usepackage{ytableau}
\ytableausetup{aligntableaux=center}

\usepackage{amsxtra}
\usepackage{eucal}

\newtheorem{rmk}{Remark}

\author{D. Korb}
\title{Order on the fixed points of the Gieseker variety with respect to the torus action}

\begin{document}

\maketitle

%\begin{center}
%\large{\textbf {Order on the fixed points of the Gieseker variety with respect to the torus action}}
%\end{center}

\section{Introduction}

In this paper I consider three families of orders on the set of $r$-multipartitions of $n$, depending on parameters from $\mathbb Q^r$. The one we are interested in is a geometric order on the set of fixed points for a torus action on the Gieseker variety. The other two serve as bounds for the geometric one, providing an almost accurate description.

Below, all three orders are defined, and it is proven that the combinatorial ones give upper and lower bounds for the geometric one. Then it is shown that in some cases all three coincide, but it is not so in general.

For references on quiver varieties and Gieseker varieties in particular, one may consult Nakajima's works \cite{Nak1} and \cite{Nak2}.

I must thank I.Loseu for proposing this problem to me and for multiple useful discussions on the subject.

%%%%%%%%%%%%%%%%%%%%%%%%%%%%%%%%%%%%%%%%%%%%%%%%%%%%%%%%%%%%%%%%%%%%%%%%%%%%%%%%%%%%
\section{Definitions}

Fix a pair of natural numbers $(n, r)$. An $r$-multipartition of $n$ is a family of partitions $\Lambda = (\lambda^1, \lambda^2, \dots, \lambda^r)$ such that $\sum |\lambda^i|=n$. These partitions will be called \textit{components} of a multipartition. A multipartition can be viewed as an array of $r$ Young diagrams, and we'll use this description to speak of boxes of a multipartition.

We introduce three different families of orders on the set of $r$-multipartitions, parametrized by vectors $(\chi_1, \chi_2, \dots, \chi_r) \in \mathbb Q^r$, two of which are combinatorial and, in fact, very similar, and the last one is geometric, coming from the orders on the set of fixed points with respect to certain torus action on the Gieseker moduli space. From this point of view, the vector $\chi$ can be considered as a character of a torus action on the framing at infinity.

All these constructions should be considered for a generic $\chi$. $\chi$ is called generic if and only if the combinatorial pre-orders described below are non-degenerate (it follows easily from the construction that this condition is the same for both). This condition is equivalent to the following: $\chi_i-\chi_j \notin \{0, 1, \dots, n-1\}$ $\forall i\neq j$ Geometrically it corresponds to the action with the finite number of fixed points.

%%%%%%%%%%%%%%%%%%%%%%%%%%%%%%%%%%%%%%%%%%
\subsection{Combinatorial orders}

Let us introduce a shifted content function $Cont_{\chi}$ on boxes of multipartitions as follows: it equals $\chi_i$ in the topmost left box of $\lambda^i$, it increases by 1 with every box when moving to the right, and decreases by 1 when moving down.

The first order is given as follows: $\Lambda \geq M$ if and only if the $n$-vector of values of $Cont_{\chi}$ on the boxes of $\Lambda$ is greater than or equal to the same on the boxes of $M$ in the following sense: there exists a bijection between the boxes such that the shifted content of the box of $\Lambda$ is greater than or equal to the shifted content of the corresponding box of $M$. Alternatively, one can order the boxes of $\Lambda$ $b_1, \dots, b_n$ and the boxes of $M$ $b_1', \dots, b_n'$ so that $Cont_{\chi}(b_i) \geq Cont_{\chi}(b_i')$ for every $i$.

The second order is somewhat more complicated. It may be considered a generalisation of the usual dominance order on the partitions defined via the procedure of dropping the boxes. Consider two multipartitions $\Lambda$ and $M$. They are called \textit{adjacent} if two sets of boxes $\Lambda \setminus M$ and $M \setminus \Lambda$ form a skew partition each (in particular, each of these sets of boxes belong to one component of multipartition), and these skew partitions coincide up to the choice of the component and translation. For adjacent $\Lambda$ and $M$ we say that $\Lambda \triangleright M$ if $\Lambda \geq M$. Completing this relation with respect to transitivity, we get a new order relation $\Lambda \triangleright M$. It is obvious from this definition that $\Lambda \triangleright M$ implies $\Lambda \geq M$.

%%%%%%%%%%%%%%%%%%%%%%%%%%%%%%%%%%%%%%%%%%
\subsection{Geometric order}

Consider a Gieseker moduli space $\mathcal M_{(r, n)}$, which we describe, following Nakajima, as a quiver variety. Namely, consider two vector spaces $V$, $W$, $\dim V=n$, $\dim W=r$, and $M_{(r, n)}=\{(B_1, B_2, i, j)| [B_1, B_2]+ij=0, Span \left<B_1^a B_2^b \operatorname{Im}(i) \right>=V\}\subset End(V) \times End(V) \times Hom(W, V) \times Hom(V, W)$, then $\mathcal M_{(r, n)}=M_{(r, n)}/GL(V)$, where $GL(V)$ acts naturally by base change in $V$. The natural projection $\pi: M_{(r, n)} \to \mathcal M_{(r, n)}$ is a principal $GL(V)$-bundle.

We consider two tori acting on $M_{(r, n)}$, and, consequently, on $\mathcal M_{(r, n)}$. First, we fix a maximal torus $\mathbb T_f$ in $GL(W)$, acting naturally. Second, we consider a two-dimensional torus $\mathbb T_h$, acting by multiplication of $B_i$ by scalars. Consider a one-dimensional subtorus $\rho: \mathbb C^* \hookrightarrow \mathbb T_f \times \mathbb T_h$ such that $\rho(t).B_1=t^kB_1$, $\rho(t).B_2=t^{-k}B_2$ and consider points fixed by $T=\mathbb T_f \times \mathbb T_h$ in $\mathcal M_{(r, n)}$. For a fixed point $p$, let $D(p)=\{x \in \mathcal M_{(r, n)}| \lim_{t \to \infty} \rho(t).x = p\}$. Now, for two fixed points $p \succeq q$ if $q \in \overline{D(p)}$. It is not an order as is, so we have to complete it with respect to transitivity. Note that $\rho$ is determined by $k$ and $\chi \in Hom(\mathbb C^*, \mathbb T_f) \simeq \mathbb Z^r$, and the order remains the same if we take $\rho^l$ instead of $\rho$, so, replacing $(\chi, k)$ with $\frac{\chi}{k}$, we obtain a family of orders parametrized by $\chi \in \mathbb Q^r$.

Now we shall describe the set of fixed points. One can find this description in \cite{Nak3}, section 3.1, but we provide an explicit construction here as we will need it later.

For every fixed point $p$, there exists a homomorphism $y: \mathbb T_f \times \mathbb T_h \to GL(V)$ and a point $P=(B_1, B_2, i, j) \in M_{(r, n)}$ with $\pi(P)=p$ such that the natural action of $\mathbb T_f \times \mathbb T_h$ on $P$ coincides with the one induced by $y$. In particular, limiting ourselves to the action of $\mathbb T_h$, we get a $\mathbb Z^2$-grading on $V$ such that $B_1$ and $B_2$ have degrees $(1, 0)$ and $(0, 1)$, respectively. Let us consider $W$ having degree $(0, 0)$ so that $i$ is of degree $(0, 0)$ too. Due to the semistability condition $Span \left<B_1^a B_2^b \operatorname{Im}(i) \right>=V$, we have $V$ living in the $\mathbb Z_{\geq 0}^2$ part of the grading. Since $[B_1, B_2]+ij=0$, the map $j$ must be of degree $(1, 1)$. As $W$ has degree $(0, 0)$ and $V_{(-1, -1)}=0$, we obtain that $j=0$ for every lift $P$ of every $\mathbb T_f \times \mathbb T_h$-fixed point $p \in \mathcal M_{(r, n)}$. It follows immediately that $[B_1, B_2]=0$.

For an $r$-multipartition $\Lambda$ of $n$ let $p_{\Lambda}$ denote a point constructed as follows: we identify $V$ with a vector space spanned by the vectors indexed by the boxes of $\Lambda$ and fix vectors $w_1, w_2, \dots, w_r \in W$, which are eigenvectors corresponding to the choice of maximal torus $\mathbb T_f \subset GL(W)$. Now we fix operators $(B_1, B_2, i)$ as follows: $i(w_l)$ equals to the vector, corresponding to the topmost left corner of $\lambda^l$ ($0$ if the latter is empty), $B_1$ acts by moving the boxes one step to the right and $B_2$ acts by moving the boxes one step downwards (again, should it lead us out of a diagram, we declare the image to be equal to $0$). It is easy to see that this point is fixed with respect to $\mathbb T_f \times \mathbb T_h$. This construction can be reversed, so that for every fixed point we obtain a uniquely defined multipartition. Namely, for every $w_l$ all nonzero vectors of the form $B_1^a B_2^b i(w_l)$ can be parametrized by the boxes of a Young diagram, so that $B_1$ and $B_2$ act as prescribed. Hence this geometric order can be considered an order on the set of multipartitions, which is again denoted $\succeq$.

For representation theoretic applications, it is important to have an order on the set of points fixed by a one-parametric subgroup. One can show that in the case of non-generic $\chi$ one has infinitely many such points while a generic subgroup fixes precisely the same points as the large torus. For this reason we consider only the case of generic $\chi$.

%%%%%%%%%%%%%%%%%%%%%%%%%%%%%%%%%%%%%%%%%%
\subsection{Examples}

For partitions, I use the following notation: $\lambda = (\lambda_1, \lambda_2, \dots, \lambda_k)$, where $\lambda_i$ denotes the length of the $i^{th}$ row of the Young diagram. Here is an example for $\lambda = (3, 3, 1)$:

$$\ydiagram{3,3,1}$$

For a multipartition the notation is as before, $\Lambda = (\lambda^1, \lambda^2, \dots, \lambda^r)$ and $\lambda^i=(\lambda^i_1, \lambda^i_2, \dots, \lambda^i_k)$. Here is an example of a multipartition $\Lambda = ((2, 1), (0), (1, 1, 1), (2))$:

$$\left(\ydiagram{2,1}, \emptyset, \ydiagram{1,1,1}, \ydiagram{2}\right)$$

For $\chi=(5, 2, \frac32, -2)$ the boxes of this multipartition have the following shifted contents:

$$\left(\ytableaushort{56,4}, \emptyset, \ytableaushort{{\frac32}, {\frac12}, {-\frac12}}, \ytableaushort{{-2}{-1}} \right)$$

Here are two examples of multipartitions adjacent to $\Lambda$: $\Lambda'=((2, 2), (0), (1, 1), (2))$ and $\lambda''=((2, 1), (1, 1), (1), (2))$:

\begin{multline*}
\Lambda'=\left(\ydiagram{2, 2}*[*(gray)]{1+0, 1+1}, \emptyset, \ydiagram{1,1}, \ydiagram{2}\right)\\
\Lambda''=\left(\ydiagram{2,1}, \ydiagram[*(gray)]{1, 1}, \ydiagram{1}, \ydiagram{2}\right)
\end{multline*}

Both $\Lambda'$ and $\Lambda''$ are greater than $\Lambda$ with respect to the combinatorial orders. Note that they are not adjacent to each other.

Let us describe the fixed point corresponding to $\Lambda$. To do this, we describe a point $P=(B_1, B_2, i, j) \in M_{(4, 8)}$ such that its equivalence class $[P] \in \mathcal M_{(r, n)}$ is a fixed point, or, equivalently, there is a map $s_{\Lambda}: \mathbb T_f \times \mathbb T_h \to GL(V)$ such that $s(t)P=t.P$. For any choice of basis ${e_1, e_2, \dots, e_8}$ in $V$, we establish the following correspondence between the vectors of the basis and the boxes of $\Lambda$:

$$\left(\ytableaushort[e_]{12,3}, \emptyset, \ytableaushort[e_]{4, 5, 6}, \ytableaushort[e_]{78} \right)$$

Our construction yields, then, that we must choose the following point:

\begin{multline*}
P=(B_1=e_1^* \otimes e_2+e_7^* \otimes e_8, B_2=e_1^* \otimes e_3 + e_4^* \otimes e_5 + e_5^* \otimes e_6,\\ i=w_1^* \otimes e_1 + w_3^* \otimes e_4 + w_4^* \otimes e_7, j=0)
\end{multline*}

Here the notation $u^*$ is used for an element of the dual basis corresponding to a basis vector $u$ and $Hom(A, B)$ is identified with $A^* \otimes B$. The corresponding action of $\mathbb T_f \times \mathbb T_h$ on $V$ is then given by the following weights:

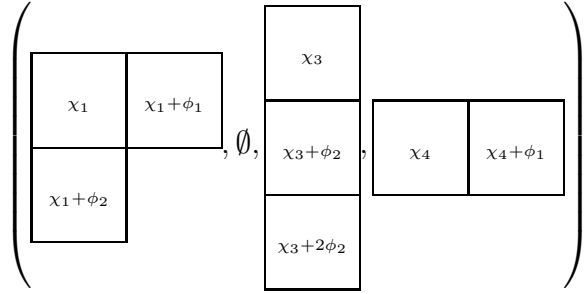
\ytableausetup{boxsize=3em}

$$\left({\tiny \ytableaushort[\scriptstyle \chi_1]{{}{+\phi_1},{+\phi_2}}, \emptyset, \ytableaushort[\scriptstyle \chi_3]{{}, {+\phi_2}, {+2\phi_2}}, \ytableaushort[\scriptstyle \chi_4]{{}{+\phi_1}} }\right)$$

Here $\chi_i$ is the character of $\mathbb T_f$ on $w_i$ while $\phi$ is the character of $\mathbb T_h$. In particular, for given a 1-parametric subgroup $\rho$ as before, the torus weight on every basis vector coincide with the shifted content of the corresponding box. This is easily seen to be the general case.

%%%%%%%%%%%%%%%%%%%%%%%%%%%%%%%%%%%%%%%%%%%%%%%%%%%%%%%%%%%%%%%%%%%%%%%%%%%%%%%%%%%%
\section{Upper and lower bounds}

Fix $\chi \in \mathbb Q^r$. In this section I intend to prove two claims:

\begin{enumerate}

\item $\Lambda \triangleright M \Rightarrow \Lambda \succeq M$.
\item $\Lambda \succeq M \Rightarrow \Lambda \geq M$.

\end{enumerate}

%%%%%%%%%%%%%%%%%%%%%%%%%%%%%%%%%%%%%%%%%%
\subsection{Upper bound}

We shall prove that $\Lambda \succeq M \Rightarrow \Lambda \geq M$.

Suppose that $\Lambda \ngeq M$. This is tautologically equivalent to the following: for every bijection between the boxes of $\Lambda$ and the boxes of $M$, there exist such a pair of corresponding boxes $b\in \Lambda$, $b'\in M$ that $Cont_{\chi}(b)<Cont_{\chi}(b')$. Let $p_{\Lambda}$ and $p_M$ denote the fixed points corresponding to $\Lambda$ and $M$, respectively.

I shall construct a certain line bundle on $\mathcal M_{(r, n)}$ and a section of this bundle, which is $0$ on $D(p_{\Lambda})$ and does not vanish in $p_M$. The existence of such a section implies $\Lambda \nsucceq M$.

Consider the tautological principal $GL_n$-bundle $\pi: M_{(r, n)} \to \mathcal M_{(r, n)}$ and the associated tautological vector bundle $\mathcal V$: the fiber at every point naturally corresponds to the space $V$ from the quiver description. We also have tautological morphisms $B_{\alpha}: \mathcal V \to \mathcal V$, $i: \mathcal O^r \to \mathcal V$. The line bundle in question is $\mathcal O (1)=\det \mathcal V$.

Every multipartition can be considered as a "recipe" of constructing $n$ vectors in V, more precisely, given a multipartition $N$ we can construct a morphism $\mathcal O^n \to \mathcal V$, which is an isomorphism in a neighbourhood of the corresponding fixed point. Namely, for a box of $N$ belonging to the partition $\nu^l$ with coordinates $(x, y)$ (from the topmost left corner which is $(0, 0)$), consider a vector $B_1^xB_2^yi(w_l)$. Our construction of the fixed points implies that these vectors form a basis in $V$ in the corresponding fixed point (and, consequently, in some neighbourhood of the latter). The formula $B_1^xB_2^yi(w_l)$, however, makes sense for every point, yielding a section $e_{(x, y)}^l:\mathcal O \to \mathcal V$ for every box and a morphism $\mathfrak E_N: \mathcal O^n \to \mathcal V$ for every multipartition. This construction is nothing but a parametrisation of basis in $V$ using the boxes of the multipartition.

Let us consider the line bundle $\det \mathcal V$ with the section $\det \mathfrak E_M$, which is nonzero in $p_M$ and zero in every other fixed point (for every $N \neq M$ there exists a box that lies in $M$ and not in $N$, hence the corresponding vector vanishes). Let us prove that for $\Lambda \ngeq M$ this section vanishes on $D(p_{\Lambda})$.

Consider the neighbourhood $\mathcal U$ of $p_{\Lambda}$, defined by the condition $\det \mathfrak E_{\Lambda} \neq 0$. In such a neighbourhood one can choose a basis in $V$ uniformly using $\mathfrak E_{\Lambda}$, that is, choose the basis in the same way as in $p_{\Lambda}$. Let us write $\mathfrak E_M$ in this basis, that is, consider the transition matrix $\mathfrak E_{\Lambda}^{-1} \circ \mathfrak E_M \in End(\mathcal O^n)$. The matrix coefficients are eigenvectors with respect to the torus action, and the weight of the coefficient corresponding to the boxes $a \in \Lambda, b \in M$ equals $Cont_{\chi}(b)-Cont_{\chi}(a)$. Note that these coefficients are equal to either $1$ or $0$ at $p_{\Lambda}$, depending on whether or not $b=a$. This means that for any point in $D(p_{\Lambda}) \cap \mathcal U$ any matrix coefficient with the weight greater than $0$ must be equal to $0$, otherwise the limit of $\rho(t).x$ at $\infty$ does not exist. Then the condition $\Lambda \ngeq M$ is equivalent to the condition that every summand of the determinant (represented as a sum over the set of all permutations) contains a matrix element of weight greater than $0$, which must consequently be equal to $0$. Hence, $\det \mathfrak E_M$ vanishes on $D(p_{\Lambda})$. \qed

%%%%%%%%%%%%%%%%%%%%%%%%%%%%%%%%%%%%%%%%%%
\subsection{Lower bound}

We shall prove that $\Lambda \triangleright M \Rightarrow \Lambda \succeq M$.

Consider two adjacent multipartitions $\Lambda \triangleright M$. We will construct a one-dimensional orbit $\mathbb C^*p' \subset \mathcal M_{(r, n)}$, such that $\lim_{t \to \infty} t.p' = p_{\Lambda}$ and $\lim_{t \to 0} t.p' = p_M$. This would imply $\Lambda \succeq M$.

The boxes in the intersection $\Lambda \cap M$ will be called \textit{fixed} and the rest will be called \textit{moved}. Moved boxes such that their top and left neighbours were also moved will be called \textit{inner} and the rest will be called \textit{border}. The sets of fixed, moved, inner and border boxes in $\Lambda$ will be denoted $\Lambda^{fix}$, $\Lambda^{mov}$, $\Lambda^{in}$ and $\Lambda^{bd}$ respectively.

\textbf{Example.}
This is an example of two adjacent multipartitions. Fixed boxes are white, inner moved boxes are black and border boxes are gray.
\ytableausetup{boxsize=1em}
$$\Lambda = \left(\dots \ydiagram{5, 5, 4, 1}*[*(gray)]{3+2, 2+1, 2+1}*[*(black)]{0, 3+2, 3+1}, \dots, \ydiagram{7, 2, 1, 1, 1}, \dots \right)$$
$$M = \left(\dots \ydiagram{3, 2, 2, 1}, \dots, \ydiagram{7, 4, 4, 3, 1}*[*(gray)]{0,2+2,1+1, 1+1}*[*(black)]{0, 0, 2+2, 2+1}, \dots \right)$$

For a multipartition $\Lambda$ consider again the choice of a basis as before, $\mathfrak E_{\Lambda}: \mathcal O^n \to \mathcal V$. The choice of basis gives us a section of the principal bundle $\pi$ in the neighbourhood of $p_{\Lambda}$: identification of $V$ with the fiber of $\mathcal V$ and writing every map in coordinates allows us to choose, for every point $q \in \mathcal U$, an element $\phi_{\Lambda}(q) \in M_{(r, n)}$ such that $\pi( \phi_{\Lambda}(q))=q$. In particular, $\phi_{\Lambda}(p_{\Lambda})$ is the same point we used to initially define the point $p_{\Lambda}$. As before, $e_{(x, y)}^l=B_1^xB_2^yi(w_l): \mathcal O \to \mathcal V$ is the basis vector corresponding to the box $(x, y)$ inside the $l^{th}$ component of the multipartition.

Let us construct the point $p'$. Consider the basis $\mathfrak E_{\Lambda}$ and let $\phi_{\Lambda}(p_{\Lambda})=(B_1, B_2, i, 0)$. The point $p'$ will have the form $\pi((B_1, B_2, i, 0)+C)$ where $C$ depends on $\Lambda$ and $M$. If the set of moved boxes coincides with one of the partitions $\mu_k$ from $M$, let $\epsilon=w_k^* \otimes e_b$ and put $C=(0, 0, \epsilon, 0)$, where $e_b$ is the  basis vector in $\mathfrak E_{\Lambda}$, corresponding to the top left moved box and $u^*$ is again an element of the dual basis. Otherwise, put $C=(\varepsilon_1, \varepsilon_2, 0, 0)$ where $\varepsilon_{1, 2}$ are computed as follows. For every box $b \in \Lambda^{bd}$, if there is a fixed box above the corresponding box in $M$, let $e_b^u$ denote the corresponding basis vector. Likewise, if there is a fixed box to the right of the corresponding box in $M$, let $e_b^r$ denote the corresponding basis vector. Put $\varepsilon_1=\sum_{b \in \Lambda^{bd}} (e_b^r)^* \otimes e_b$, $\varepsilon_2=\sum_{b \in \Lambda^{bd}} (e_b^u)^* \otimes e_b$.

This construction means that we add several ``missing'' matrix units which are present in the corresponding maps for $p_M$ to $B_1$, $B_2$ and $i$. In particular, replacing $\Lambda$ with $M$ in this construction and vice versa we obtain the same point in $\mathcal M_{(r, n)}$.

Note that $\rho(t)(p')=\pi((B_1, B_2, i, 0)+t^{-d}C)$, where $d$ is the difference in shifted contents between any of the moved boxes in $\Lambda$ and the corresponding box in $M$. By the assumption, $d>0$, and exchanging $\Lambda$ and $M$ we would get $d<0$. In conjunction with the previous observation this yields $\lim_{t \to \infty} t.p' = p_{\Lambda}$ and $\lim_{t \to 0} t.p' = p_M$. \qed

\begin{rmk}
Note that for non-generic values of $\chi$ one could choose two different adjacent partitions with the same shifted content. In this case the points $\pi((B_1, B_2, i, 0)+tC)$ for $t \in \mathbb C^*$ will be fixed with respect to $\rho$, so the geometric order does not make much sense in this situation.
\end{rmk}

%%%%%%%%%%%%%%%%%%%%%%%%%%%%%%%%%%%%%%%%%%%%%%%%%%%%%%%%%%%%%%%%%%%%%%%%%%%%%%%%%%%%
\section{Comparison of the combinatorial orders}

Now I will show that the combinatorial orders coincide inside a certain ``asymptotic'' chamber, and provide an example where these two orders differ.

%%%%%%%%%%%%%%%%%%%%%%%%%%%%%%%%%%%%%%%%%%
\subsection{Asymptotic case}

Consider a weight vector in the asymptotic chamber $A=\{\chi \in \mathbb Q^r|\chi_i-\chi_{i+1}>n-1\}$. Let $\Lambda$ and $M$ be two multipartitions, such that $\Lambda > M$. Let us show that one can move one box in $\Lambda$ so that for the resulting multipartition $\Lambda'$ holds $\Lambda \triangleright \Lambda' \geq M$. This would imply that $\Lambda \geq M$ if and only if $\Lambda \triangleright M$.

Let us call a box of a Young diagram $removable$, if the remaining part of the diagram is again a Young diagram.

Set $\Lambda_{nk+i}$ denote $\lambda^{k+1}_i$, considering every empty row of a partition having length $0$. In the case of the asymptotical chamber as above, one can see that the order $\geq$ can be described as follows: $\Lambda \geq M$ if and only if for any $l$ such that $1 \leq l \leq nr$, we have $\sum_{i=1}^l (\Lambda_i-M_i) \geq 0$.

Note that if there exists $k$ such that $|\lambda^i|=|\mu^i| \forall i<k$ and $|\lambda^k|>|\mu^k|$, one can drop the lowermost removable box from $\lambda^k$ into the first row of $\lambda^{k+1}$, and the inequality $\Lambda \geq M$ will be preserved. Note that this is the place where the choice of the chamber is important. Otherwise, if there is no such $k$ then $|\lambda^i|=|\mu^i|$ and $\lambda^i \geq \mu^i$ $\forall i$ with respect to the usual dominance order on partitions. It is well known, then, that the two descriptions of this order coincide.

%%%%%%%%%%%%%%%%%%%%%%%%%%%%%%%%%%%%%%%%%%
\subsection{Example}

Consider $n=6$, $r=4$ and a weight vector $\chi=(0, c, 2+a, 2+b)$ for  arbitrary $0<a<b<c<1$. Let us consider two multipartitions: $\Lambda = ((0), (3), (1^2), (1))$ and $M = ((3), (0), (1), (1^2))$. Note that $\Lambda \geq M$. Note also that, for every adjacent $\Lambda' \triangleleft \Lambda$, $\Lambda'$ is either smaller than or not comparable to $M$ with respect to the ordering $\geq$. Hence, $\Lambda \geq M$ does not always imply $\Lambda \triangleright M$.

\ytableausetup{boxsize=2.5em}

$$Cont_{\chi}(\Lambda):\left({\tiny \emptyset, \ytableaushort{{c}{1+c}{2+c}}, \ytableaushort{{2+a},{1+a}}, \ytableaushort{{2+b}}} \right)$$

$$Cont_{\chi}(M):\left({\tiny \ytableaushort{012}, \emptyset, \ytableaushort{{2+a}}, \ytableaushort{{2+b},{1+b}}} \right)$$

Note that this example is in some sense minimal, that is, while this example is not unique, there are no such examples for smaller $n$ or $r$.

%%%%%%%%%%%%%%%%%%%%%%%%%%%%%%%%%%%%%%%%%%%%%%%%%%%%%%%%%%%%%%%%%%%%%%%%%%%%%%%%%%%%


\begin{thebibliography}{99}

\bibitem[Na1]{Nak1} H. Nakajima, \textit{Instantons on ALE spaces, quiver varieties, and Kac-Moody algebras}, Duke, vol. 76, no.2 (November 1994)

\bibitem[Na2]{Nak2} H. Nakajima, \textit{Lectures on Hilbert scheme of points on surfaces}, AMS (1999)

\bibitem[NY]{Nak3} H. Nakajima, K. Yoshioka, \textit{Lectures on instanton counting}, CRM Proceedings and Lecture Notes, Proceedings of "Workshop on algebraic structures and moduli spaces" (July 14-20, 2003)

\end{thebibliography}
\end{document}